\documentclass[12pt]{amsart}
\usepackage{amsmath}
\usepackage{amsfonts}
\usepackage{amsthm}
\usepackage{amssymb}
\usepackage{amscd}
\usepackage{epic}
\usepackage{eepic}

\newcommand{\printname}[1]
  {\smash{\makebox[0pt]{pace{-2.0in}\raisebox{8pt}{\tiny #1}}}}
\newtheorem {Theorem}   {Theorem}
\newtheorem {Lemma} {Lemma} [section]
\newtheorem {Proposition}[Lemma]{Proposition}

\theoremstyle{definition}
\newtheorem{Definition}[Lemma]{Definition}
\theoremstyle{remark}
\newtheorem{Remark}[Lemma]{Remark}

\newtheorem {Corollary}[Lemma]{Corollary}
\newcommand {\cal}[1]  {{\mathcal{#1}}}
\newcommand{\ssminus}{{\smallsetminus}}

\newcommand{\nab}[1][]{\ensuremath{\mathrm{\nabla}{#1}}}

\newcommand{\be}{\begin{equation}}
\newcommand{\ee}{\end{equation}}

\newcommand{\lb}[1]{\label{#1}}

\newcommand{\bet}{\beta}
\newcommand{\al}{\alpha}
\newcommand{\om}{\omega}

\newcommand{\we}{\wedge}

\newcommand{\lam}{\lambda}

\newcommand{\Ref}[1]{(\ref{#1})}
\newcommand{\fr}{\frac}

\newcommand{\ra}{\rightarrow}

\newcommand{\ri}{\mathrm{r}}
\newcommand{\Lie}{\pounds}
\newcommand{\Lap}{\Delta}

\newcommand{\wht}{\widehat}
\newcommand{\sk}{s.k.r.p.}
\def\t{\tau}

\newcommand{\sol}{Ricci soliton}
\newcommand{\sols}{Ricci solitons}

\newcommand{\Mt}{M_\t}

\begin{document}
\title{Special K\"ahler-Ricci potentials and Ricci solitons}


\author{Gideon Maschler}

\address{Department of Mathematics and Computer Science, Emory University,
Atlanta, Georgia 30322, U.S.A.}
\email{gm@mathcs.emory.edu}



\begin{abstract}
On a manifold of dimension at least six, let $(g,\t)$  be a pair consisting of
a K\"ahler metric $g$ which is locally K\"ahler irreducible, and a nonconstant 
smooth function $\t$. Off the zero set of $\t$, if the metric $\wht{g}=g/\t^2$ 
is a gradient \sol\ which has soliton function $1/\t$, we show that 
$\wht{g}$ is K\"ahler with respect to another complex structure, and locally of a 
type first described by Koiso. Moreover, $\t$ is a special K\"ahler-Ricci potential, 
a notion defined in earlier works of Derdzinski and Maschler. The result extends
to dimension four with additional assumptions. We also discuss a {\em  
Ricci-Hessian equation}, which is a generalization of the soliton equation, and 
observe that the set of pairs $(g,\t)$ satisfying a Ricci-Hessian equation is 
invariant, in a suitable sense, under the map $(g,\t)\ra (\wht{g},1/\t)$.
\end{abstract}

\maketitle

\setcounter{Theorem}{0}
\renewcommand{\theTheorem}{\Alph{Theorem}}
\renewcommand{\theequation}{\arabic{section}.\arabic{equation}}

\section{Introduction}
\setcounter{equation}{0}

In this paper we study  pairs $(g,\t)$ on a manifold $M$ of dimension larger than two, 
where $g$ is a Riemannian metric and $\t$ is a smooth nonconstant function. In this
context, an important role will be played by the map 
$(g,\t)\ra (\wht{g},\wht{\t})=(g/\t^2,1/\t)$ on the set of all such pairs.
On $M\ssminus\t^{-1}(0)$, this map is a well-defined involution. We will call 
$(\wht{g},\wht{\t})$ the {\em associated pair} of $(g,\t)$.

We say that a pair satisfies a {\em Ricci-Hessian equation} if
\be\lb{ric-hessian} \al\nab d\t +\ri=\gamma\, g\end{equation}
holds for the Hessian of $\t$, the Ricci tensor $\ri$ of $g$, and 
some $C^\infty$ coefficient functions $\al$ and $\gamma$. If $\al$ and $\gamma$
are constant, the pair, or sometimes just the metric, is called a gradient \sol.

Before stating our main result, we note a closely related fact: the set of pairs 
satisfying a Ricci-Hessian equation is invariant under our involution. The latter is
well-defined once the domain of each pair is further restricted (see \S\ref{dual}).
In this setting we call it the duality map.
 
A part of our main result may be stated informally as follows. Consider two subsets 
of the set of all pairs on $M$: those for which the metric is K\"ahler (and locally 
irreducible in a suitable sense, see below), and pairs which are gradient \sols.  
Assuming a restriction on the dimension of the manifold, if the involution maps an 
element of the first subset to an element of the second one, the latter element lies 
in the intersection of the two subsets. More precisely,
\begin{Theorem}\lb{main}
Let $M$ be a manifold of dimension at least six, and $(g,\t)$ a pair as above, 
with $g$ a K\"ahler metric. Suppose $g$ is not a local product of 
K\"ahler metrics in any neighborhood of some point of $M$. 
If the associated pair $(\wht{g},\wht{\t})$ is a gradient \sol, then, 
on $M\ssminus \t^{-1}(0)$, the metric $\wht{g}$ is K\"ahler.
\end{Theorem}

The complex structures giving the K\"ahler structures of $g$ and $\wht{g}$ are 
oppositely oriented. Also, with an extra assumption, the result extends to real 
dimension four.

The proof, in fact, yields much more information on both pairs.
First, $\t$ is a special K\"ahler-Ricci potential. This notion (Definition 
\ref{skrp-def}) was first defined in \cite{local}, and implies that $\t$ is a Killing 
potential, and that \Ref{ric-hessian} holds in some (generically nonempty) open set 
of the manifold. 
Second, the K\"ahler-\sol\ $(\wht{g},\wht{\t})$ is locally of a type first described 
by Koiso in \cite{ko} (see also \cite{ca}). 

Although our result is of a local character, one should note that there exist compact  
manifolds, specifically toric Fano ones, which admit K\"ahler-\sols\ \cite{wz}, 
most of which are not of the form found by Koiso.

In the following we describe a few related problems of a broader scope. The involution 
above is defined in part via a conformal change, and one can ask whether, starting with 
a K\"ahler metric $g$, one can find a metric $g/\t^2$, for some function $\t$ as above, 
which is a \sol. The case where $g/\t^2$ is Einstein, was the subject of the study
of \cite{local,skrp,global}, where local and global classifications were given, 
and, in all even dimensions larger than four, $\t$ turns out to be, in fact, a 
special K\"ahler-Ricci potential. In dimension four this need not be the case, and 
different compact examples where recently shown to exist in \cite{clw}. 

For the  case of \sols, no such general results are
known, even if one {\em assumes} that $\t$ is a special K\"ahler-Ricci potential. 
Our result can be considered a first step in an attempt to answer this question. 
For more on this topic, see the material surrounding Proposition \ref{inv-one}.

In various talks, G. Tian has asked whether there exist compact non-K\"ahler \sols\
in dimension four. Extending the question to all dimensions, one may answer it
affirmatively via constructions involving products. Ignoring these fairly simple
examples, one can try to produce such a \sol\ in the conformal class of a K\"ahler 
metric (at least in dimension four, it is not too difficult to see that there
can be no more than two K\"ahler metrics in a given conformal class). Our result can
be regarded as implying that, in a special case, such an attempt will fail.

Finally, recall the result of Schur, stating that if $\ri=\phi\, g$ for some function
$\phi$, then, except in dimension two, $\phi$ must be constant. A similar principle
holds for K\"ahler-\sols, for the coefficient of $g$, and one may ask whether it
holds for any \sol, or, equivalently, whether there exist pairs for which 
\Ref{ric-hessian} holds with the coefficient of $\nab d\t$ constant, but not the 
coefficient of $g$. In unpublished work, A. Derdzinski has shown that such pairs do 
exist on compact manifolds. In fact, on these manifolds there are families of pairs 
$(g,\t)$, with $g$ K\"ahler, for which the associated pairs $(\wht{g},\wht{\t})$ are 
each of this  type, and in fact are obtained by deforming one of the Einstein metrics in 
\cite{global}.

This paper is arranged as follows. Riemannian preliminaries on duality and \sols\ 
appear in \S\ref{rds}. Consequences of the K\"ahler condition for \Ref{ric-hessian}, 
along with a review of the basic properties of metrics with a special 
K\"ahler-Ricci potential, are given in \S\ref{kahl}. Ordinary differential equations 
associated with a large class of such metrics are studied in \S\ref{he}, especially in 
relation to the assumption that the associated pair forms a \sol, and we give an 
analysis of their solution set. After recalling the geometric structure of a K\"ahler 
metric admitting a special K\"ahler-Ricci potential, and presenting
a duality result in this context in \S\ref{geo-dua}, we summarize our results in 
\S\ref{sol-light} by proving Theorem \ref{detail}, from which Theorem \ref{main} easily 
follows. Our conventions throughout closely follow \cite{local}.


\setcounter{Theorem}{0}
\renewcommand{\theTheorem}{\thesection.\arabic{Theorem}}

\section{Ricci-Hessian equations, duality and \sols}\lb{rds}
\setcounter{equation}{0}
 
\subsection{conformal changes}
Let $(M,g)$ be a Riemannian manifold of dimension $n$, and 
$\t:M\to {\mathbb{R}}\,$ a nonconstant $\,C^\infty$ function. We  write metrics
conformally related to $g$ in the form $\wht{g}=g/\t^2$, and we set $Q=g(\nab\t,\nab\t)$.
We always consider the metric $\wht{g}$ on its domain of definition, i.e. the set 
$M\ssminus\t^{-1}(0)$. With respect to $\wht{g}$, the Hessian of a given $C^2$ function 
$f$ on $M\ssminus\,\t^{-1}(0)$ is given by

\begin{equation}\lb{hess}
\begin{array}{l}
\wht\nabla df\,=\,\nabla df\,+\,\t^{-1}[2\,d\t\odot df
-\,g(\nab\t,\nab f)g],
\end{array}
\end{equation}
where $d\t\odot df=(d\t\otimes df+df\otimes d\t)/2$. 
We will be concerned primarily with the case where $df\we d\t=0$, i.e., at points where 
$d\t\neq 0$, $f$ is given locally as a composition $f=H\circ\t$. In this case, 
\Ref{hess} becomes $\wht\nabla df\,=f'\,\nabla d\t\,
+\left(f''+2\t^{-1}f'\right)\,d\t\otimes d\t-f'\,\t^{-1}\,Q\,g$, with 
$'$ denoting differentiation with respect to $\t$. 
For the particular choice of $f=\t^{-1}$, this expression simplifies: 
\be\lb{observe}
\wht\nabla d\t^{-1}\,=-\t^{-2}(\,\nabla d\t\,-\,\t^{-1}\,Q\,g),\quad\mathrm{\ if\ } 
\wht{g}=g/\t^2
\end{equation}

We conclude by recording the conformal change expression relating the
Ricci tensors of $g$ and $\wht{g}$, with $\Lap$ denoting the Laplace operator:
\be\lb{Ricci}
\wht{\ri}\,=\,\ri\,+\,(n-2)\,\t^{-1}\nabla d\t\,+\,
\left[\t^{-1}\Delta\t\,-\,(n-1)\,\t^{-2}Q\right]g.\\
\end{equation}

\subsection{Ricci-Hessian equations and duality}\lb{dual}
With $M$, $g$, $\t$ and other notations as above, we say that the pair $(g,\t)$ 
satisfies a {\em Ricci-Hessian equation} on $M$ (or often just on an open set of $M$), 
if \Ref{ric-hessian} holds there. We record this equation more explicitly as
\be\lb{ric-hes}\al\nab d\t +\ri=\gamma\, g,\quad \mathrm{with\ } \t 
\mathrm{\ nonconstant},
\end{equation}
where $\nab d\t$  and $\ri$ are as above, and $\al$, $\gamma$ are $C^\infty$ 
coefficient functions. What we will call duality may be regarded informally
as an involution on the space of pairs satisfying \Ref{ric-hes}:
\begin{Proposition}\lb{duality}
Let $M$ have dimension $n>3$, and suppose a pair $(g,\t)$ as above
satisfies a Ricci-Hessian equation \Ref{ric-hes} on $M$. Then the 
pair $(\wht{g},\wht{\t})=(g/\t^2,1/\t)$ also satisfies a 
Ricci-Hessian equation $\wht{\al}\wht{\nab} d\t +\wht{\ri}=\wht{\gamma}\, \wht{g}$,
on $M\ssminus\t^{-1}(0)$, with coefficients 
\be\lb{dual-coef}\wht{\al}=(n-2)\t-\t^2\al,\quad\wht{\gamma}=
\gamma\t^2-(1+\al\t)Q+\t\Lap\t.
\end{equation}
\end{Proposition}
In fact, letting $\bet$ denote the coefficient of $g$ in \Ref{Ricci}, one has, 
by \Ref{observe} and \Ref{Ricci},  
\begin{eqnarray*}
\wht{\al}\,\wht{\nabla} d\wht{\t}+\,\wht{\ri}&=&
((n-2)\t-\t^2\al)(-\t^{-2}\,\nabla d\t\,+\t^{-3}Qg)+\,\ri\,
+(n-2)\t^{-1}\nab d\t+\bet g\\
&=&\al\,\nabla d\t\,+\,\ri\,
+(\wht{\al}\t^{-3}Q+\bet)g=(\gamma+\wht{\al}\t^{-3}Q+\bet)\t^2\wht{g},
\end{eqnarray*}
and one easily checks that the last expression is $\wht{\gamma}\,\wht{g}$.

\begin{Remark}
As mentioned in the introduction, the pair $(\wht{g},\wht{\t})$ will be called the
associated pair. It is not necessarily defined on all of $M$, 
and it also does not determine the coefficients $\al$, $\gamma$ uniquely at every
point of $M$. Hence, to consider $(g,\t)\ra (\wht{g},\wht{\t})$ as an involution on 
the set of metrics satisfying \Ref{ric-hes} for some smooth 
coefficients $\al$ and $\gamma$, one has to restrict consideration to
the complement in $M$ of $\t^{-1}(0)$ and the closed set of points where 
$\nab d\t$ is a multiple of $g$. In some cases one may also consider coefficients
with isolated singularities, in which case those singularities also must be excluded

To verify the involutive property, one easily checks that $\wht{\wht{\al}}=\al$, while 
$\wht{\wht{\gamma}}=\gamma$ follows from this, as $g$, $\t$  and $\alpha$ determine 
$\gamma$ uniquely. One can also check the last relation directly, using the following 
formulas for the two functions $\wht{Q}=
\wht{g}({\wht{\nab}\wht{\t}},{\wht{\nab}\wht{\t}})$ and $\wht{\Lap}\wht{\t}$:
$$\wht{Q}=\t^{-2}Q,\qquad  \wht{\Lap}\wht{\t}=n\t^{-1}Q-\Lap\t.$$  
\end{Remark}
\begin{Remark}\lb{affine}
For any pair $(g,\t)$ satisfying \Ref{ric-hes}, one can produce another such 
pair by an affine change in $\t$ (a nonconstant one). If this affine change involves 
only a change by an additive constant, it leaves equation \Ref{ric-hes} invariant.
This freedom induces, of course, a freedom in the choice of $\wht{\t}$ , which will be 
exploited in Proposition \ref{skrp-dual}.
\end{Remark}


\subsection{Ricci solitons}
A {\it \sol\/} \cite{hamilton-tr} is a Riemannian manifold $\,(M,\wht{g})\,$ 
such that $\,\Lie_v \wht{g}+\,\wht{\ri}\,=e\,\wht{g}$ for some constant $e$ and 
$\,C^\infty$ vector field $\,v\,$ on $\,M$. Here $\,\Lie_v$ is the Lie derivative and 
$\,\wht{\ri}\,$ denotes the Ricci tensor of $\wht{g}$.
We will only be interested in {\em gradient} \sols, in which $M$
admits a $\,C^\infty$ function $\,f:M\to {\mathbb{R}}\,$ with
\begin{equation}\label{nfr}
\wht{\nabla} df\,+\,\wht{\ri}\,=\,e\, \wht{g}\quad\mathrm{for\ a\ constant\ }
e.
\end{equation}
We will call $f$ the {\em soliton function}.
By a result of Perelman \cite[Remark~3.2]{p}, every compact 
\sol $\,(M,\wht{g})\,$ is a gradient \sol. Recall also that a metric
is Einstein if its Ricci tensor is a multiple of it.

Thus a gradient \sol\ is nothing but a pair $(\wht{g},f)$ satisfying a 
Ricci-Hessian equation with constant coefficients.
Using \Ref{observe} and \Ref{Ricci}, or, more naturally, 
the duality of  Proposition \ref{duality} (slightly modified to
allow $\wht{\t}$ to be multiplied by a constant), we have
\begin{Proposition}\lb{sol-opt}
Let $(M,g)$ be a Riemannian manifold of dimension $n>2$
and $\t$ a nonconstant $\,C^\infty$ function. The \sol\ equation 
$\wht{\nab} d\,(b\t^{-1})\,+\,\wht{\ri}\,=\,e\, \wht{g}$,
with $b$ a constant, 
holds for $\wht{g}=\t^{-2}g$ on $M\ssminus\t^{-1}(0)$, if and only if
$g$ satisfies a Ricci-Hessian equation \Ref{ric-hes} with coefficients
\be\lb{no-symmetric}
\al=(n-2)\,\t^{-1}-b\,\t^{-2}, \quad
\gamma=e\,\t^{-2}-\t^{-1}\Delta\t
+\left((n-1)\,\t^{-2}-b\,\t^{-3}\right)Q.
\end{equation}
\end{Proposition}

\begin{Remark}\lb{non-Ein}
The introduction of the constant $b$ is not, strictly speaking, necessary
for later development, but serves to compare with the conformally Einstein
case, which occurs when $b=0$: relations \Ref{no-symmetric} with $b=0$
are implied by \cite[(6.1) and (6.2)]{local}, which hold in that case.
It follows from this that an Einstein metric cannot also satisfy
a \sol\ equation \Ref{nfr} with the soliton function a {\em nonzero} multiple 
of $\tau^{-1}$. Note that for other nonconstant soliton functions, this is 
possible.
\end{Remark}

\begin{Remark}
Here and in Proposition \ref{inv-one} we briefly consider 
the conformal change equation for $g$, which yields a 
gradient \sol\ $\wht{g}$ with an arbitrary soliton function $f$. 
An analogous calculation using equation \Ref{Ricci} yields
\begin{eqnarray}\lb{gen-f}
\ri&+&(n-2)\t^{-1}\nab d\t\,
+\nabla df+2\t^{-1}\,d\t\odot df\nonumber\\
&=&\left[e\t^{-2}+(n-1)\t^{-2}Q-
\t^{-1}\Delta\t+\t^{-1}g(\nab\t,\nab f)\right]g.
\end{eqnarray}
If $df\we d\t=0$, this gives
\begin{eqnarray}\lb{f(tau)}
\ri&+&\left(f'+(n-2)\,\t^{-1}\right)\,\nab d\t\,
+\left(f''+2\t^{-1}f'\right)\,d\t\otimes d\t\nonumber\\
&=&\left[e\,\t^{-2}-
\t^{-1}\Delta\t+\left((n-1)\t^{-2}+\t^{-1}f'\right)Q\right]g.
\end{eqnarray}
We can conclude from this that a particular choice for $f$ will 
eliminate the Hessian term. Namely, setting
$f=-(n-2)\log|\t|$, the metric  $\wht{g}$ is a \sol\
metric precisely when
$$
\ri-\left(n-2\right)\t^{-2}\,d\t\otimes d\t
=\left[\t^{-2}\left(e+Q\right)-
\t^{-1}\Delta\t\right]g\,.
$$
However, while this equation looks quite simple, it implies that
$g$ cannot be K\"ahler (unless $n=2$ or $\t$ is constant). 
This is another reason, apart from duality considerations, why we will 
focus on the case of a soliton function proportional
to $\t^{-1}$.
\end{Remark}

\section{Ricci-Hessian equations and special K\"ahler-Ricci potentials}\lb{kahl}
\setcounter{equation}{0}
\subsection{The Ricci-Hessian equation and \sols} 
Let $(M,J)$ be a complex manifold, with $J$ the associated almost complex structure. 
Suppose $g$ is a K\"ahler metric on $M$, i.e a Riemannian metric for which $J$ 
is parallel. Let $(g,\t)$ be a pair satisfying the Ricci-Hessian equation \Ref{ric-hes} 
on $M$. The K\"ahler property implies that both $g$ and $\ri$ are Hermitian, hence so is 
$\nab d\t$ on the support of $\al$. In many of our applications, this support will be 
a dense set in $M$. This property of $\nab d\t$ is equivalent to the statement that 
$\t$ is a Killing potential, i.e. a $C^\infty$ function for which $J\nab\t$ is a 
Killing vector field (cf. 
\cite[Lemma 5.2]{local}).  

In the K\"ahler case, if $g/\t^2$ is a \sol, certain restrictions on the soliton
function force it to be proportional to, or at least affine in $\t^{-1}$.
\begin{Proposition}\lb{inv-one}
Let $(M,g)$ be a K\"ahler manifold with a Killing potential $\t$, and $\wht{g}=g/\t^2$ 
a \sol\ with a $\t$-dependent soliton function $f$. Then $f$ is an affine function in 
$\t^{-1}$. 
\end{Proposition}
\begin{proof}
As \Ref{f(tau)} holds under our assumptions, and $d\t\otimes d\t$ is the only term
in it that is not Hermitian, its coefficient $f''+2\t^{-1}f'$ must vanish, implying 
the conclusion.
\end{proof}
\begin{Remark}
If the Killing assumption above is replaced by \Ref{ric-hes} for some 
$\t$, $\al$ not identically zero and $\gamma$, the conclusion still follows
on the support of $\al$. If one then drops the $\t$-dependency assumption on the soliton 
function $f$, all that \Ref{ric-hes} implies, in combination with \Ref{gen-f}, is that 
$\nab df+2\t^{-1}d\t\odot df=(2\t^2)^{-1}\Lie_{(\t^2\nabla f)}\, g$ is Hermitian on the 
support of $\al$.
\end{Remark}

We will be especially interested in the case where $\al$ and $\gamma$ in \Ref{ric-hes}
are functions of $\t$. We note that, this always holds for $\al$ in \Ref{no-symmetric}, 
while it holds for $\gamma$ there if both $d\t\we d\Lap\t=0$ and $d\t\we dQ=0$.
One may attempt to weaken these assumptions using methods akin to those of
\cite[(6.5) and Proposition 6.4]{local}. We choose to follow here the quicker approach
of \cite[\S1.4]{ham1}, which, however, works only for $m>2$.
\begin{Proposition}\lb{coeff}
If \Ref{ric-hes} holds for a K\"ahler metric of complex dimension $m>2$, 
and $d\al\we d\t=0$, then $d\gamma\we d\t=0$.
\end{Proposition}
\begin{proof}
Composing \Ref{ric-hes} with $J$ and applying $d$ to the result gives
$d\al\we d(\imath_{\nab\t}\om/2)=(d\al/d\t)\,d\t\we d(\imath_{\nab\t}\om/2)=
d\gamma\we\om$, using \cite[(5.3)]{local} (here $\om$ is the K\"ahler form of $g$). 
Exterior multiplication with $d\t$ gives $d\t\we d\gamma\we\om=0$, and the result 
follows because the operation $\we\om$ is injective on $2$-forms for $m>2$.
\end{proof}

Thus the coefficients of the Ricci-Hessian equation $\Ref{no-symmetric}$ will be 
functions of $\t$, provided $(M,g,\t)$ is K\"ahler of dimension $m>2$, and $g/\t^2$ is 
a \sol\ with soliton function proportional to $\t^{-1}$. 

\subsection{Special K\"ahler-Ricci potentials}\lb{skrps}

Below, we denote by $\Mt$ the complement, in a manifold $M$, of the critical set of a 
smooth function $\t$. For a Killing potential on a K\"ahler manifold, $\Mt$ is 
open and dense in $M$.
\begin{Definition}\lb{skrp-def}\cite{local}
A nonconstant Killing potential $\t$ on a K\"ahler manifold $(M,J,g)$ is 
called a {\em special K\"ahler-Ricci potential} if, on the set $\Mt$,
all non-zero tangent vectors orthogonal to $\nab\t$ and $J\nab\t$ are 
eigenvectors of both $\nab d\t$ and $\ri$.
\end{Definition}
We will call a metric admitting a special K\"ahler-Ricci potential a {\em \sk\ metric}, 
and occasionally will declare $(g,\t)$ to be a {\em \sk\ pair}. Among the more important 
characteristics of such a metric is the existence of a Ricci-Hessian equation.
More precisely
\begin{Proposition}\lb{skrp-hess}\cite[Corollary 9.2, Remarks 7.1 and 7.4]{local}
Let $(M,g)$ be a K\"ahler manifold of complex dimension $m\geq 2$. If 
\Ref{ric-hes} holds for some $C^\infty$ functions $\al$, $\gamma$ and (nonconstant) 
$\t$, with $d\al\we d\t=0$, $d\gamma\we d\t=0$ and $\al\,d\al\ne 0$ everywhere in $\Mt$, 
then $\t$ is a special K\"ahler-Ricci potential. 
Conversely, if $(M,g)$ admits a K\"ahler-Ricci potential $\t$, then \Ref{ric-hes} 
holds on an open subset of $\Mt$, namely away from points where $\nab d\t$ 
is a multiple of $g$.
\end{Proposition}
\begin{Remark}
In \cite{local}, we have actually written the Ricci-Hessian equation in the form  
$\nab d\t +\chi\ri=\sigma\, g$. Note that the domains of the coefficient functions
may vary as one switches between these two forms. In general, any statement involving  
the Ricci-Hessian equation of a \sk\ metric refers to the largest 
domain on which \Ref{ric-hes} holds. Moreover, this change results in a slightly 
different statement of the first part of Proposition \ref{skrp-hess}, while to get the
second part (and its proof), one need only to switch $r$  with $\nab d\t$ in 
\cite[second paragraph of Remark 7.4]{local}.
\end{Remark}
\begin{Corollary}\lb{sol-skrp-2}
If $(M,g)$ is K\"ahler, of complex dimension $m>2$, and $\wht{g}=g/\t^2$ is
a \sol, with soliton function $b\t^{-1}$, where $b$ is a constant, then
$\t$ is a special K\"ahler-Ricci potential.
\end{Corollary}
\begin{proof}
We combine Propositions \ref{sol-opt}, \ref{coeff} and \ref{skrp-hess},
except that for $\al$ in \Ref{no-symmetric}, $\al\, d\al=0$ on the set
where $\t=(n-2)/b$ and on the set where $\t=2b/(n-2)$, so that $\t$
is a special K\"ahler-Ricci potential, and hence a Killing potential
away from these sets. By \cite[Lemma 5.2]{local}, $\nab d\t$ is
Hermitian away from these sets, yet it is also clearly Hermitian
in the interior of each of these sets, so that by continuity, it is
Hermitian on all of $\Mt$. Again using \cite[Lemma 5.2]{local}, this means 
that $\nab\t$ is holomorphic on $\Mt$, which implies that the interiors 
of the above mentioned two sets are empty. As the \sk\ condition is
defined by equalities, continuity now implies that $\t$ satisfies it 
on all of $\Mt$. 
\end{proof}

By \cite[Definition $7.2$, Remark $7.3$]{local},
the \sk\ condition on $(g,\t)$ is equivalent to the existence, 
on $\Mt$, of an orthogonal decomposition 
$TM={\cal{V}}\oplus {\cal{H}}$, with
${\cal{V}}=\mathrm{span}(\nab\t, J\nab\t)$, along with
four smooth functions  $\phi$, $\psi$, $\lam$, $\mu$
which are pointwise eigenvalues for either $\nab d\t$ or $\ri$, 
i.e., they satisfy
\begin{equation}\lb{skrp}
\begin{array}{rclrcl}
 \nab d\t|_{\cal{H}}\ &=&\phi\, g|_{\cal{H}}, \quad 
&\nab d\t|_{\cal{V}}\ &=&\psi\, g|_{\cal{V}},\\
 \ri|_{\cal{H}}\ &=&\lam\, g|_{\cal{H}}, \quad
&\ri|_{\cal{V}}\ &=&\mu\, g|_{\cal{V}}.
\end{array}
\end{equation}
This decomposition is also $\ri$- and $\nab d\t$-orthogonal.

\begin{Remark}\lb{nontriv}
By \cite[Lemma 12.5]{local}, $\phi$ either vanishes identically on 
$\Mt$, or never vanishes there. In the former case,
$g$ is reducible to a local product of K\"ahler metrics near any point
(see \cite[Corollary 13.2]{local} and \cite[Remark 16.4]{skrp}). 
In the latter case, we call $g$ a {\em nontrivial} \sk\ metric. 
\end{Remark}

\begin{Remark}\lb{c-const}
For a nontrivial \sk\ metric, consider $c=\t-Q/(2\phi)$,  
with $Q=g(\nab d\t,\nab d\t)$, and $\kappa=\mathrm{sgn}(\phi)(\Lap\t+\lam\,Q/\phi)$, 
regarded as functions $\Mt\ra \mathbb{R}$. By \cite[Lemma 12.5]{local}, 
$c$ is constant on $\Mt$, and will be called the {\em \sk\ constant}.
In any complex dimension $m\geq 2$, we will call a nontrivial \sk\ metric 
{\em standard} if $\kappa$ is constant (and also use ``standard \sk\ pair"
as a designation for $(g,\t)$). According to \cite[\S27, using (10.1) 
and Lemma 11.1]{local}, constancy of $\kappa$ holds if $m>2$, so that the 
designation ``standard" involves an extra assumption as compared with ``nontrivial" 
only when $m=2$. The geometric meaning of $\kappa$ will be recalled in \S\ref{geo}.
\end{Remark}

\begin{Remark}\lb{chi-relate}
Using Proposition \ref{skrp-hess}, for any \sk\ metric, \Ref{ric-hes} holds
on points of $\Mt$ for which $\phi\neq\psi$. On this subset, we
regard \Ref{ric-hes} as an equation of operators, and equate eigenvalues to obtain 
$\,\al\phi+\lam=\gamma=\al\psi+\mu$, so that 
\be\lb{alpha}\,\lam-\mu=(\psi-\phi)\al.\end{equation}
According to \cite[Lemma 11.1a]{local}, $Q$, $\Lap\t$, $\phi$, $\psi$ and $\mu$ 
are locally $C^\infty$ functions of $\t$ on $\Mt$. If $g$ is a standard
\sk\ metric, $\lam$ is also such a function, as one concludes from the equation defining 
$\kappa$. Hence, by \Ref{alpha}, the same holds for $\al$ on its domain of definition.
\end{Remark}


\section{Associated differential equations}\label{he}
\setcounter{equation}{0}
\subsection{The \sk\ differential equation} 

A number of ordinary differential equations are associated with 
nontrivial \sk\ metrics. Special cases of these were given in
\cite{local}. They are derived below from the Ricci-Hessian
equation \Ref{ric-hes}, i.e.
$$\al\nabla d\t\,+\,\ri\,=\,\gamma\, g.$$

In the next proposition, $\al$ will be as in \Ref{ric-hes}, $\phi$ as in \Ref{skrp},
$c$ and $\kappa$ as in Remark \ref{c-const} and a prime denotes the derivative operator 
$d/d\t$.
\begin{Proposition}\lb{diff-eq}
Let $(g,\t)$ be a \sk\ pair with $g$ nontrivial, on a manifold $M$ of 
complex dimension $m$. Then, the equation
\begin{equation}\lb{mek}
(\t-c)^2\phi''+\,(\t-c)[ m-(\t-c)\al]\phi'-\,m\phi\,=\,-\mathrm{sgn}(\phi)\kappa/2.
\end{equation}
holds at points of $\Mt$ for which  $\phi'(\t)$ is nonzero. 
If $g$ is standard, \Ref{mek} is an ordinary differential equation, which, upon 
differentiation and division by $\t-c$, takes the homogeneous form
\begin{equation}\label{tcp}
(\t-c)\phi'''\,
=\,\,[(\t-c)\al-m-2]\,\phi''\,+\,[(\t-c)\al'+2\al]\,\phi'.
\end{equation}
\end{Proposition}
A special case of equation \Ref{tcp} was important in \cite{global}, but is given here
mainly for the sake of completeness. We will only be using equation \Ref{mek}.
\begin{proof}
By Remark \ref{chi-relate}, on $\Mt$, each of
$Q$, $\phi$, $\psi$, $\Lap\t$ and $\mu$ is locally a function
of $\t$. In fact, we have
\begin{equation}\label{qet}
\begin{array}{rrclrlcl}\arraycolsep5pt
\mathrm{a)}&\psi\,&=&\,\phi\,+\,(\t-c)\phi',\phantom{_{j_j}}\quad
&\mathrm{b)}&\psi'&=&\,2\phi'+\,(\t-c)\phi'',\\
\mathrm{c)}&\Lap\t&=&\,2m\phi\,+\,2(\t-c)\phi', \quad
&\mathrm{d)}&\mu\,&=&\,-(m+1)\phi'-\,(\t-c)\phi''.
\end{array}
\end{equation}
Namely, \cite[Lemma 11.1(b)]{local} gives
$\,2\psi=Q'$, which yields (\ref{qet}.a) (and hence (\ref{qet}.b)),
since $\,Q=2(\t-c)\phi\,$ due to the definition of $\,c$. Next,
(\ref{qet}.c) is immediate from (\ref{qet}.a), as
$\,\Lap\t=\mathrm{tr}_g\nabla d\t=2\psi+2(m-1)\phi$. Finally,
$\,2\mu=-(\Lap\t)'$ by \cite[Lemma 11.1(b)]{local}, and so,
differentiating (\ref{qet}.c), we obtain (\ref{qet}.d).

Next, by Remark \ref{chi-relate}, equation \Ref{ric-hes} holds on points of 
$\Mt$ for which $\phi\neq\psi$. Since, using (\ref{qet}.a), 
the latter inequality holds when $(\t-c)\phi'\neq 0$, and $\t\neq c$ on 
$\Mt$ (as $Q = 2 (\tau - c) \phi$  and  $Q > 0$  on $\Mt$), 
we see that this set consists exactly of the points of $\Mt$ for which 
$\phi'(\t)$ is nonzero.

As on this subset of $\Mt$, \Ref{ric-hes} holds, so does
\Ref{alpha}, which along with $Q = 2 (\tau - c) \phi$  and the definitions of 
$\,\kappa\,$ and $\,c\,$  easily yields $\,\mathrm{sgn}(\phi)\kappa/2=\Lap\t/2+
(\t-c)\lam=\Lap\t/2+(\t-c)[\mu+(\psi-\phi)\al]$. Replacing $\,\mu,\psi\,$ and
$\,\Lap\t$ with the expressions provided by (\ref{qet}), we get (\ref{mek}).
If $g$ is standard, $\kappa$ is constant, so $\al$ is a function of $\t$ by Remark 
\ref{chi-relate}. Hence  equation \Ref{mek} is an ordinary differential equation, and
\Ref{tcp} then follows as described in the body of the proposition.
\end{proof}
\begin{Remark}\lb{cnvrs}
A converse statement to this result can be made, where \Ref{mek} implies \Ref{ric-hes}
for a standard \sk\ metric, under the following extra assumptions.

Let $\,\phi\,$ be {\it globally\/} a function of $\,\t$, in
the sense that it is the composite of $\,\t\,$ with some $\,C^\infty$
function $\,I'\to {\mathbb{R}}\,$ on the image interval $\,I'=\t(\Mt)$.
(That $I'$ is indeed an interval is known, see \cite[\S10 and \S11]{skrp}.)
Assuming $\,\phi'$, as a function of $\,\t$, is nonzero at all points of a
dense subset of $\,I'$, and (\ref{mek}) holds on $\,I'$ for
a $\,C^\infty$ function $\,\al:I'\to {\mathbb{R}}$, it follows that (\ref{ric-hes}) is
satisfied on $\,\Mt$ by  $\al=\al(\t)$ and some $\,\gamma$.

In fact, the assumption involving $\phi'$ means, as we have seen in the proof above,
that \Ref{ric-hes} holds on a dense subset of $\Mt$,
with some $\,\al\,$ that must coincide with the one above:
they both satisfy (\ref{mek}) with the same $\,\phi$ on a dense subset of
$\,I'$, and hence everywhere in $\,I'$.
\end{Remark}

\subsection{The differential equations in relation to \sols}\lb{dif-sol}
Let $(g,\t)$ be a standard \sk\ pair on a manifold $M$ of complex dimension $m$. 
Suppose $g/\t^2$ is a Ricci soliton with soliton function $b\t^{-1}$, 
where $b$ is a constant. By Proposition \ref{sol-opt}, equation \Ref{ric-hes} 
holds on $M\ssminus\t^{-1}(0)$, with 
\be\lb{al-sol}
\al=(2(m-1)\,\t-b)/\t^2.
\end{equation}
Hence, in this case, the differential equations appearing in Proposition \ref{diff-eq} 
take the form 
\be\lb{sol1}
\t^2(\t-c)^2\phi''+(\t-c)\left[m\t^2-(\t-c)(2(m-1)\t-b)\right]\phi'-m\t^2\phi
=-\mathrm{sgn}\,(\phi)\kappa\,\t^2/2
\end{equation}
and
$$\t^3(\t-c)\phi'''=[(m-4)\t^3-(2(m-1)c+b)\t^2+bc\t]\phi''+
[2(m-1)\t(\t+c)-2bc]\phi'.$$
These equations hold for $\t$ values corresponding to points of $M_\t$ on which 
$\phi'(\t)$ is nonzero.

Another ordinary differential equation is obtained on the same set as follows. 
The term $\gamma$ in \Ref{ric-hes} is given, by Remark \ref{chi-relate}, as
$\gamma=\al\psi+\mu$, and we substitute for $\psi$ and $\mu$ their respective 
expressions (\ref{qet}.a) and (\ref{qet}.d), to obtain 
$$\gamma=\al\phi+\left(\al(\t-c)-(m+1)\right)\phi'-(\t-c)\phi''.$$
In the case at hand, $\gamma$ also has an expression derived from the last
term of \Ref{no-symmetric}, in which we replace $Q$ by $2(\t-c)\phi$, and 
$\Lap\t$ by (\ref{qet}.c). Equating the two expressions, and replacing $\al$ by
\Ref{al-sol}, we get after rearranging terms and multiplying by $\t$ that
\begin{eqnarray}\lb{sol3}
-\t^3(\t-c)\phi''&+&[(2m\t-b)\t(\t-c)-\t^3(m+1)]\phi'\nonumber\\
&+&[2(2m-1)\t^2-b\t+2(\t-c)(b-(2m-1)\t)]\phi=e\t.
\end{eqnarray}
The fact shown shortly, that \Ref{sol3} is not, in general, a consequence of \Ref{sol1}, 
is the main local difference between the case where the \sk\ metric is conformal to a 
non-Einstein \sol\ of the type we are considering, with $b\neq 0$, and the one where 
it is conformal to an Einstein metric ($b=0$). The latter was the object of study of 
\cite{local,skrp,global}.

\subsection{Solutions of the system \{\Ref{sol1}, \Ref{sol3}\}}

To examine the solutions of the system \{\Ref{sol1}, \Ref{sol3}\}, 
we note the following
\begin{Lemma}\lb{nonrat}
Let $\{\phi'+p\phi=q,\, A\phi''+B\phi'+C\phi=D\}$ be a system of ordinary differential
equations in the variable $\t$, with coefficients $p$, $q$, $A$, $B$, $C$ and $D$ 
that are rational functions. Then, on any nonempty interval admitting a solution
$\phi$, either
\be\lb{rat-rel} A(p^2-p')-Bp+C=0
\end{equation} holds identically, or
\be\lb{sol-doub}\phi=\left(D-A(q'-pq)-Bq\right)/\left(A(p^2-p')-Bp+C\right).
\end{equation} holds away from the (isolated) singularities of the right hand side.
\end{Lemma}
\begin{proof}
Let $\phi$ be a solution on an interval as above. We have 
$\phi'=q-p\phi$, so that $\phi''=q'-p'\phi-p\phi'=q'-p'\phi-p(q-p\phi)=
(p^2-p')\phi+q'-pq$. Substituting this in the second equation, while collecting 
terms involving $\phi$, gives
$$\left(A(p^2-p')-Bp+C\right)\phi+A(q'-pq)+Bq=D,$$
from which the result follows at once.
\end{proof}

To apply this lemma to \{\Ref{sol1}, \Ref{sol3}\}, we multiply
\Ref{sol3} by $(\t-c)/\t$ and add it to \Ref{sol1}, giving a first order
equation which, after multiplying by $\t$, simplifies to 
\begin{eqnarray}\lb{1st-ord}
\t^2(\t-c)(\t-2c)\phi'\!\!&+&\!\!\left[-m\t^3+(2(2m-1)c+b)\t^2-(2(2m-1)c^2+3bc)\t+
2bc^2\right]\!\!\phi\nonumber\\
&=&-\mathrm{sgn}(\phi)\kappa\t^3/2+e\t(\t-c).
\end{eqnarray}
We will be applying Lemma \ref{nonrat} to the system \{\Ref{1st-ord}, \Ref{sol1}\}
(modifying \Ref{1st-ord} appropriately). This system has, 
a solution set identical to that of \{\Ref{sol1}, \Ref{sol3}\}
(certainly on intervals not containing $0$, $c$ and $2c$, and by a continuity
argument, on any interval). To compute \Ref{rat-rel} in this case, note that $p$, 
given as a ratio of the coefficient of $\phi$ to that of $\phi'$ in \Ref{1st-ord},  
has a partial fraction decomposition of the form 
$p=m/(\t-c)-1/(\t-2c)+b/\t^2-(2m-1)/\t$. Similarly, 
$q=\mathrm{sgn}(\phi)\kappa/(2(\t-c))
+(e-2\mathrm{sgn}(\phi)\kappa c)/(2c(\t-2c))-e/(2\t c)$. Using $A=\t^2(\t-c)^2$, 
$B=(\t-c)[m\t^2-(\t-c)(2(m-1)\t-b)]$,  $C=-m\t^2$ and 
$D=-\mathrm{sgn}(\phi)k\t^2/2$, two long but direct computations gives
\begin{equation}\lb{rati}
\begin{array}{lcl}\arraycolsep5pt
D-A(q'-pq)-Bq &=& 0,\\
A(p^2-p')-Bp+C&=&-2bc(\t-c)^2/(\t(\t-2c)).
\end{array}
\end{equation} 
This immediately gives
\begin{Proposition}\lb{bc}
Suppose $bc\ne 0$. Then the system \{\Ref{sol1}, \Ref{sol3}\} has no nonzero solutions
on any nonempty open interval.
\end{Proposition}
\begin{proof}
Assume $bc \ne 0$. Then the second of equations \Ref{rati} does not vanish identically 
on the given interval. Hence Lemma \ref{nonrat} implies that any solution is the ratio
of the two equations in \Ref{rati}, away from the point $c$. This ratio is the zero
function. By continuity, \{\Ref{1st-ord}, \Ref{sol1}\}, as well as 
\{\Ref{sol1}, \Ref{sol3}\}, admit no nonzero solutions on the given interval.
\end{proof}

\subsection{Solutions for the case $c=0$}
If $c=0$, equations \{\Ref{sol1},\,\Ref{sol3}\} take the form, 
\begin{equation}\lb{ode13}
\begin{array}{rcrclcl}\arraycolsep5pt
\t^4\phi''&+&\left[(2-m)\t^3+b\t^2\right]\phi'&-&m\t^2\phi
&=&-\mathrm{sgn}\,(\phi)\kappa\,\t^2/2,\\
-\t^4\phi''&+&[((m-1)\t^3-b\t^2]\phi'&+&b\t\phi&=&e\t,
\end{array}
\end{equation}
with $m$ a positive integer, and $b\ne 0$.
As special solutions, one can take $\mathrm{sgn}(\kappa)\kappa/(2m)$ for the first, and 
$e/b$ for the second. A basis of solutions to each associated homogeneous equation is
given by $\{\t^m\exp(b/\t),\sum b^{m-l}\t^l/(m-l)!\}$, where the sum ranges over
$l=0\ldots m-1$ for the first, and $l=1\ldots m-1$ for the second. Thus, the general
solution to the system has the form  
\be\lb{phi-sol}
\phi=A+B\sum_{l=1}^{m-1}\frac{b^{m-l}}{(m-l)!}\t^l+
C\t^m\exp({\,b/\t})
\end{equation}
for arbitrary constants $A$, $B$ and $C$ (where $A$ represents the sum of an arbitrary 
multiple of $b^m/m!$ with $\mathrm{sgn}(\kappa)\kappa/(2m)+e/b$).


\section{Geometry and duality for \sk\ metrics}\lb{geo-dua}
\setcounter{equation}{0}

\subsection{Local Geometry of \sk\ metrics}\lb{geo}
We recall here the main case in the geometric classification of \sk\ metrics.
Let $\pi:(L,\langle\cdot,\cdot\rangle)\ra (N,h)$ be a hermitian holomorphic line bundle 
over a K\"ahler-Einstein manifold of complex dimension $m-1$. Assume that the 
curvature of $\langle\cdot,\cdot\rangle$ is a multiple of the K\"ahler form of $h$. 
Note that, if $N$ is compact and $h$ is not Ricci flat, this implies that $L$ is smoothly 
isomorphic to a rational power of the anti-canonical bundle of $N$.

Consider, on $L\ssminus N$ (the total space of $L$ excluding the zero section), 
the metric $g$ given by
\be\lb{metric-form}
g|_{\cal{H}}=2|\t-c|\,\pi^*h, \quad g|_{\cal{V}}=\fr {Q(\t)}{(ar)^2}\,
\mathrm{Re}\,\langle\cdot,\cdot\rangle, 
\end{equation}
where\\
-- $\cal{V},\cal{H}$ are the vertical/horizontal distributions of $L$, respectively,
the latter determined via the Chern connection of $\langle\cdot,\cdot\rangle$,\\
-- $c,a\neq 0$ are constants,\\
-- $r$ is the norm induced by $\langle\cdot,\cdot\rangle$,\\ 
-- $\t$ is a function on $L\ssminus N$, obtained by composing on $r$ another
function, denoted via abuse of notation by $\t(r)$, and obtained as follows:
one fixes an open interval $I$ and a positive $C^\infty$ function $Q(\t)$ on $I$, 
solves the differential equation $(a/Q)\,d\t=d(\log r)$ to obtain a diffeomorphism 
$r(\t):I\ra (0,\infty)$, and defines $\t(r)$  as the inverse of this diffeomorphism. 

The pair $(g,\t)$, with $\t=\t(r)$, is a \sk\ pair (see \cite[\S8 and
\S16]{local}, as well as \cite[\S4]{skrp}), and $|\nab\t|^2_g=Q(\t(r))$. 
If $g$ is nontrivial, the connection on $L$ will not be flat.
The constant $\kappa$ of Remark \ref{c-const} is the Einstein constant of $h$, 
so that if $g$ is nontrivial, it is in fact standard (for an arbitrary \sk\ 
metric, $h$ need not be Einstein if $m=2$). For $g$ standard, or merely nontrivial,
the \sk\ constant $c$ (see again Remark \ref{c-const}) coincides with $c$ of 
\Ref{metric-form}.

Conversely, for any standard nontrivial \sk\ metric $(M,J, g,\t)$, any point in 
$\Mt$ has a neighborhood biholomorphically isometric to an open set in 
some triple $(L\ssminus N,g,\t(r))$ as above (this is a special case of 
\cite[Theorem 18.1]{local}). This biholomorphic isometry identifies 
$\mathrm{span}\,(\nab\t, J\nab\t)$ and its orthogonal complement, with $\cal{V}$ 
and, respectively, $\cal{H}$. Moreover, one can extend $(g,\t(r))$ to all of $L$, and 
then the biholomorphic isometry can also be defined on neighborhoods of points in 
$M\ssminus \Mt$ \cite[Remark 16.4]{skrp}.

\subsection{Duality for \sk\ metrics}\lb{sk-dual}

By Proposition \ref{skrp-hess}, a \sk\ pair $(g,\t)$ satisfies a Ricci-Hessian equation 
\Ref{ric-hes}, on points of the noncritical set $M_\t$ in which $\nab d\t$ is not
a mutliple of $g$. On this set (with $\t^{-1}(0)$ excluded), the involution of 
\S\ref{dual} yields a new pair $(\wht{g},\wht{\t})$, which also satisfies a 
Ricci-Hessian equation. In general, not much can be said about $\wht{g}$. However, a 
special case of the affine change mentioned in Remark \ref{affine} involves changing 
$\t$ by an additive constant. This produces a new Killing potential $t$, with  
$(g,t)$ a \sk\ pair very closely related to $(g,\t)$. If the additive 
constant is chosen to be minus the \sk\ constant $c$, applying the involution to 
$(g,t)$ yields metrics which are K\"ahler with respect to an oppositely 
oriented complex structure. In fact, they are even \sk\ metrics. We provide a proof of 
this in the following proposition, for the sake of completeness. Similar less detailed
statements appear in \cite[Remark 28.4]{global} and \cite[end of \S5.5 and \S5.6]{ham1}.
\begin{Proposition}\lb{skrp-dual}
Let  $g$ be a standard \sk\ metric on $(M,J)$, with Killing potential $\t$
and corresponding \sk\ constant $c$. If $t=\t-c$, then the associated pair 
$(\wht{g},\wht{t})=(g/t^2, 1/t)$ is a standard \sk\ pair on $M\ssminus\t^{-1}(c)$.
\end{Proposition}
In fact, the proof will imply that the metric $\wht{g}$ is K\"ahler with respect 
to the complex structure $\bar{J}$ given by $\bar{J}|_{\cal{H}}=J|_{\cal{H}}$, 
$\bar{J}|_{\cal{V}}=-J|_{\cal{V}}$, where ${\cal{H}}$ is the orthogonal complement 
of ${\cal{V}}=\mathrm{span}\,(\nab\t, J\nab\t)$.  This structure, defined
on $\Mt$, extends uniquely to $M$ (see \cite[Remark 16.4]{skrp}), and the 
corresponding extension of $\wht{g}$ (see end of \S\ref{geo}) to $M\ssminus\t^{-1}(0)$
is still K\"ahler with respect to it.
\begin{proof}
By the classification of \sk\ metrics, it is enough to consider $M$ as a subset of
the model line bundle $L$ of \S\ref{geo}. For simplicity, we take 
$M=L\ssminus N$. On $L$, the complex structure $\bar{J}$ defines the complex
conjugate bundle structure, which we denote $\bar{L}$. We will show that
the metric $\wht{g}$ is a \sk\ metric, by constructing it explicitly as in 
\S\ref{geo}, but on the line bundle $\bar{L}$. This line bundle is smoothly 
isomorphic to the dual bundle $L^*$, and hence the construction will transfer
to a holomorphic line bundle, which is one of the requirements for the data 
used in \S\ref{geo}. The proof that such structures are K\"ahler is indicated
in \cite[\S16]{local} (or, quite efficiently, via the methods in \cite{ham1}).

The metric $\wht{g}$ is obtained from the model metric $g$ as follows: first replace 
$\langle\cdot,\cdot\rangle$, $a$, $\t$ and $I$, respectively,  with 
the complex conjugate fiber metric $\overline{\langle\cdot,\cdot\rangle}$, 
the constant $\wht{a}=-a$, the function $\wht{t}= 1/(\t-c)$ and the
open interval $\wht{I}$ which is the image of the decreasing diffeomorphism 
$I\ni\t\ra\wht{t}\in\wht{I}$. We then replace $c$ by $\wht{c}=0$, and have 
$Q$ replaced with a function $\wht{Q}$ which is a solution to the equation 
$a\,\wht{Q}/Q=\wht{a}\,d \wht{t}/d\t$.  Finally, using these new data, along 
with $h$, $r$ and $\cal{H}$, one defines a new  standard \sk\ metric exactly as in 
\Ref{metric-form}. Note that the definition of $\wht{Q}$ guarantees that the 
required relation $(\wht{a}/\wht{Q})\,d\,\wht{t}=d(\log r)$ holds, and positivity 
of $\wht{Q}$ follows from its defining equation together with the fact that 
$\wht{t}(\t)$ is decreasing. To conclude that this standard \sk\ metric 
is indeed $\wht{g}=g/(\t-c)^2$, one computes its two factors to be
$2|\wht{t}-\wht{c}|=2/|\t-c|$ and $\wht{Q}(\wht{t})/(\wht{a}\,r)^2=
-[Q(\t)/(a\,r)^2]\,d\wht{t}/d\t=Q(\t)/[(a\,r)^2(\t-c)^2]$.
\end{proof}
\begin{Remark}\lb{phi-its-hat}
In the case  $c=\wht{c}$, i.e. $c=0$ we have $t=1/\t$, so that
$(\wht{g},t)=(\wht{g},\wht{\t})$. Then, by Remark \ref{c-const}, one has
$Q=2\t\phi$, and similarly for $\wht{Q}$. Hence $\wht{\phi}/\phi=\t\wht{Q}/(tQ)=
-(\t/t)\,d t/d\t=-(\t/(1/\t))\cdot (-1/\t^2)=1$, i.e. $\wht{\phi}(t)=\phi(\t)$. The 
same conclusion can be reached without the use of the geometric description of 
\sk\ metrics in \S\ref{geo}, by restricting \Ref{observe} to ${\cal{H}}$ and using 
\Ref{skrp} and Remark \ref{c-const}. 
\end{Remark}
\begin{Remark}
Still assuming $c=0$, and using all the above conventions, suppose one fixes all
the data defining $g$ in \Ref{metric-form}, except for $Q=2\t\phi$, which varies only 
by changing $\phi(\t)$ in the solution space of equation \Ref{mek}. If, in these
circumstances,  for each such solution, $\phi'$ satisfies the requirement in
Remark \ref{cnvrs}, it follows that \Ref{ric-hes}, and in particular, $\al=\al(\t)$
does not vary for all these metrics. As they all share the same associated equation 
\Ref{mek}, the corresponding dual metrics $\wht{g}$ also share their own 
associated equation \Ref{mek}, written with $t=\wht{\t}$ and $\wht{\al}$, the latter 
determined as in \Ref{dual-coef}. Since the solution space determines the coefficients 
of a linear differential equation, the result $\wht{\phi}(t)=\phi(\t)$ now implies that 
\Ref{mek} for $(\wht{g},\wht{\t})$ is obtained from \Ref{mek} of $(g,\t)$ simply by the
change of variable $\t\ra \wht{\t} = 1/\t$.
\end{Remark}


\section{Proof of the Theorem A}\lb{sol-light}
\setcounter{equation}{0}

\begin{Theorem}\lb{detail}
Given a standard \sk\ pair $(g,\t)$, if $\wht{g}=g/\t^2$ is a non-Einstein \sol\ with 
soliton function a multiple of $\t^{-1}$, then $\wht{g}$ is K\"ahler, and locally
of the type given by Koiso in \cite{ko} (or Cao in \cite{ca}).
\end{Theorem}
\begin{proof}
In fact, as the pair $(g,\t)$ is standard, the associated function $\phi$ cannot be 
identically zero (see remark \ref{nontriv}). But then, as a function on the image of 
$\t$, the function $\phi$ solves the system 
\{\Ref{sol1}, \Ref{sol3}\}  only if $bc=0$ (here $b\t^{-1}$ denotes, as in 
\S\ref{dif-sol}, the soliton function). As $\wht{g}$ is a non-Einstein \sol, 
$b\ne 0$ (Remark \ref{non-Ein}). Hence the \sk\ constant $c$ is zero. This implies, by Proposition \ref{skrp-dual} and the paragraph past it, that 
the soliton $\wht{g}$ is K\"ahler on $M\ssminus\t^{-1}(0)$, with respect to a complex 
structure oppositely oriented to that with respect to which $g$ is K\"ahler. By Remark 
\ref{phi-its-hat}, $\wht{\phi}(\wht{\t})=\phi(\t)$, so that, using \Ref{phi-sol} and 
the definition of $c$ in Remark \ref{c-const},
$$\wht{Q}=2\wht{\t}\wht{\phi}=\fr 2{\t}\left[
A+B\sum_{l=1}^{m-1}\frac{b^{m-l}}{(m-l)!}\fr 1{\t^l}+C\fr1{\t^m}\exp({\,b\t})\right],$$
for some constants $A,B,C$. It is  known (cf. \cite[\S2]{ham4}) that a metric 
$\wht{g}$ with the characteristics given in \S\ref{geo}, and such an expression for
$\wht{Q}$, is (locally) of the form found by Koiso.
\end{proof}
We end with the\vspace{.1in}

\noindent
{\em Proof of Theorem A.}\ 
Let $M$ be of dimension at least six, with $(g,\t)$ a pair for which $g$ is
K\"ahler. If the associated pair $(\wht{g},\wht{\t})$ is a \sol, then, by Corollary
\ref{sol-skrp-2}, $(g,\t)$ is a \sk\ pair. (This will also hold in dimension four
if $Q$ and $\Lap\t$ are $\t$-dependent, see the paragraph before Proposition 
\ref{coeff}.) The non-reducibilty assumption on $g$ implies, in these dimensions,
that it is a standard \sk\ metric. As the metric $\wht{g}$ cannot be Einstein by 
Remark \ref{non-Ein}, Theorem \ref{detail} implies that $\wht{g}$ is K\"ahler 
(and locally of the type given by Koiso).



\vspace{.3in}
\begin{center}
\textbf{Acknowledgements}
\end{center}
The author thanks A. Derdzinski for his encouragement, early
participation, reading of preliminary versions, and the numerous suggestions he made. 
These greatly improved both the style and the mathematical content of this paper, 
particularly with regard to the notion of duality.



\end{document}